\documentclass[english]{amsart}
\usepackage[T1]{fontenc}
\usepackage[latin1]{inputenc}
\usepackage{amsthm}
\usepackage{amsbsy}
\usepackage{graphicx}

\makeatletter

\numberwithin{equation}{section}
\numberwithin{figure}{section}
 \theoremstyle{definition}
 \newtheorem*{defn*}{\protect\definitionname}
\theoremstyle{plain}
\newtheorem{thm}{\protect\theoremname}
  \theoremstyle{definition}
  \newtheorem{defn}[thm]{\protect\definitionname}
  \theoremstyle{remark}
  \newtheorem*{rem*}{\protect\remarkname}

\makeatother

\usepackage{babel}
  \providecommand{\definitionname}{Definition}
  \providecommand{\remarkname}{Remark}
\providecommand{\theoremname}{Theorem}

\begin{document}

\title{FAST-SLOW VECTOR FIELDS OF REACTION-DIFFUSION SYSTEMS}

\author{V.Bykov{$^{*}$}, Y.Cherkinsky{$^{**}$}, V.Gol'dshtein{$^{**}$}, N.Krapivnik{$^{**}$},
U.Maas{$^{*}$},}

\address{{$^{*}$}Institute for Technical Thermodynamics, Karlsruhe University,
Karlsruhe, Germany}

\address{{$^{**}$}Department of Mathematics, Ben-Gurion University, Beer-Sheva,
Israel}

\begin{abstract}
A geometrically invariant concept of fast-slow vector fields perturbed
by transport terms (describing molecular diffusion)
is proposed in this paper. It is an extension of our concept of singularly
perturbed vector fields to reaction-diffusion systems.
This paper is motivated by an algorithm of reaction-diffusion manifolds (REDIM).
It can be considered as its theoretical justification extending it from a practical
algorithm to a robust computational method.
Fast-slow vector fields can be represented locally as "singularly perturbed
systems of PDE". The paper focuses on development of the decomposition
to a fast and slow subsystems. It is demonstrated that
transport terms can be neglected (under reasonable physical assumptions)
for the fast subsystem. 
A simple practical application example of the proposed algorithm for numerical treatment
of reaction-diffusion systems is demonstrated.
\end{abstract}

\maketitle

\section{Introduction}

The decomposition of complex dynamical systems into simpler subsystems using
different rates of changes (multiple time scales) for different subsystems
is common in physical and engineering models \cite{MP92}-\cite{RF01}.
The main difficulty in applications is a "hidden", implicit form of the decomposition
of the system evolving at different time scales. Namely, there is no explicit
representation of the system in relatively fast and slow subsystems available. 

A formal mathematical basis to cope with this problem is based on the notion
of singularly perturbed vector fields \cite{BGG2006}. Let us briefly introduce
some key ideas of the singularly perturbed vector fields (SPVFs) that can be
used to treat the problem of the decomposition in general.

Roughly speaking a singularly perturbed vector field (SPVF) $F(z,\varepsilon)$
is a vector field defined in a domain $G$ of Euclidian space $R^{n}$
that depends on a small parameter $\varepsilon\geq0$ such that for
any point $z$, $F(z,0)$ belongs to an a priori fixed fast subspace
$M_{f}(z)$ of smaller dimension - $\dim M_{f}(z)<n$. Moreover, the
dimension of $M_{f}(z)$ does not depend on the choice of the point
$z$. Thus, in this case the vector field $F(z,\varepsilon)$ can be decomposed
into a fast sub-field that belongs to the fast subspace $M_{f}(z)$
and its complement representing a slow sub-field. Of course this is not
a formal description, which is more sophisticated. Additionally, if $M_{f}(z)$
does not depend on $x$ then the vector field $F(z,\varepsilon)$
represents (by definition) a linearly decomposed singularly perturbed
vector field. Accordingly, the notion of the linearly decomposed singularly
perturbed vector field is a geometrical analog of a singularly perturbed
system.

A formal concept (a theory of SPVFs) can be useful for practical applications
if it is supported by an identification algorithm for these fast sub-fields \cite{BGM2008}.
In a number of previous papers an algorithm for linearly
decomposed singularly perturbed vector fields \cite{BM2009} has been constructed.
This algorithm is based on a global linear interpolation procedure for an original vector
field that we call a Global Quasi-Linearization (GQL) (see e.g. \cite{BGM2008, BM2009}). 

The theory of singularly perturbed vector fields is a coordinate
free version of singularly perturbed systems for (ODEs) \cite{GS92} - \cite{ML68}.
It cannot be used in the original form for study the influence of transport processes of
reaction-diffusion systems.
Thus, the main formal object of our current study should be modified as
$$\mathrm{F}(z,x,\varepsilon):=\Phi(z,\varepsilon)+L(z,x,\varepsilon).$$
It combines a singularly perturbed vector field $\Phi(z,\varepsilon)$ (reaction
term) and a linear operator typically of second order (diffusion term). Here $x$
belongs to a set $V_{s}$ in Euclidian space $R^{s}, s\leq3$. Typically it is a
segment $[0,L]$ or closed parallelogram. Additionally, fast reaction
terms are assumed to be much faster than corresponding transport processes, that
leads to a formal assumption $\lim_{\varepsilon\to0}L(z,x,\varepsilon)=0$,
i.e. $\mathrm{F}(z,x,0)=\Phi(z,0)$. This makes the extension
of the theoretical framework developed in \cite{BGG2006, BGM2008} straightforward.

\section{General formal Notion of Fast-Slow Vector Fields}

In this section, the main formal framework of singularly perturbed
vector fields \cite{BGG2006} to fast-slow vector fields of reaction-diffusion
systems is adopted. 

As in previous, a standard definition of vector bundles and use vector/fiber
bundles as a formal substitute for so-called nonlinear coordinate
systems.

\begin{defn*}
A vector bundle $\xi$ over a connected manifold $N\subset R^{m}$
consists of a set $E\subset R^{m}$ (the total set), a smooth map
$p:E\rightarrow N$ (the projection) which is onto, and each fiber
$F_{x}^{\xi}=p^{-1}(x)$ is a finite dimensional affine subspace.
These objects are required to satisfy the following condition: for
each $x\in N$, there is a neighborhood $U$ of $x$ in $N$, an integer
$k$ and a diffeomorphism $\varphi:p^{-1}(U)\rightarrow U\times R^{k}$
such that on each fiber $\varphi$ is an isomorphism of vector spaces. 
\end{defn*}
Note that, all fibers have to be of the same dimension $k$.
\begin{defn}
Call a domain $V\subset R^{n}$ a structured domain (or a domain structured
by a vector bundle) if there exists a vector bundle $\xi$ and a diffeomorphism
$\psi:V\rightarrow U$ onto an open subset $U\subset E$, where $E$
is the total set of $\xi$. 
\end{defn}
Fix a parametric family of smooth fast-slow vector fields 
$\boldsymbol{F}(z,x,\delta):=\Phi(z,\delta)+L(z,x,\delta)$
defined in a domain $V\subset R^{n}$ for any $0<\delta<\delta_{0}$.
Here $\delta_{0}$ is a fixed positive number and $\delta$ is a small
positive parameter (an explicit form of small parameter in the system
is needed at least at the initial stage); $\Phi(z,\delta)$ is a singularly
perturbed vector field, $L(z,x,\delta)$ is a linear differential
operator such that $\lim_{\delta\to0}L(z,x,\delta)=0$ 

A corresponding system of PDE's is then can be cast in the form
\begin{equation}
\frac{\partial z}{\partial t}=\mathbf{F}(z,x,\delta)=\Phi(z,\delta)+L(z,x,\delta).
\label{eq:main}
\end{equation}

\begin{defn}
Suppose that $V$ is a domain structured by a vector bundle $\xi$
and a diffeomorphism $\psi$. For any point $z\in G$ call $M_{z}:=\psi^{-1}(p^{-1}(\psi(z)\cap U)$
a fast manifold associated with the point $z$. Call the set of all
fast manifolds $M_{z}$ a family of fast manifolds of $V$. 

By construction any point $z\in G$ belongs to only one fast manifold.
If $z\neq z_{1}$ either $M_{z}\cap M_{z_{1}}=\emptyset$ or $M_{z}=M_{z_{1}}$.
The dimension of any manifold $M_{z}$ remains the same. Denote this dimension
by $n_{f}$ and call it the fast dimension of $G$. 

A family of fast manifolds $M_{z}$ is linear if there exists a linear
subspace $L_{f}$ of $R^{n}$ such that $M_{z}=\{z\}+L_{f}$ for any
$z\subset V$. 
\end{defn}
Call $L_{f}$ a fast subspace in this case.

This is a simplest possible "linear" situation. By using a corresponding
linear coordinate transformation of variables it is possible to move $L_{f}$ to a coordinate
subspace, which leads to the standard SPS (see e.g. \cite{F79}).

Denote by $TM_{z}$ a tangent space to $M_{z}$ at the point $z$.
\begin{defn}
A parametric family $\Phi(z,\delta):V\rightarrow R^{n}$ of vector fields
defined in a domain $V$ structured by a vector bundle $\xi$ and
a diffeomorphism $\psi$ is an asymptotic singularly perturbed vector
field if $\lim_{\delta\rightarrow0}\Phi(z,\delta)\in TM_{z}$ for
any $z\in V$ and the structure of the domain $G$ is minimal for
the vector field $\Phi(z,\delta):G\rightarrow R^{n}$ in the following
sense.

There is no a proper vector subbundle $\xi_{1}$ of the vector bundle
$\xi$ such that $\Phi(z,\delta):V\rightarrow R^{n}$ is \textit{an
asymptotic singularly perturbed vector field} in a domain $V$ structured
by the vector subbundle $\xi_{1}$ and the same diffeomorphism $\psi$.
\end{defn}
\begin{rem*}
This property of minimality means that it is not possible to reduce
the dimension of fast manifolds $\{M_{z}\}$ using sub-bundles. 
\end{rem*}
From this point outwards, without loss of generality, a family of fast 
manifolds $\{M_{z}\}$ associated with a singularly perturbed vector field $\Phi(z,\delta)$
is supposed to be minimal.

For a linear family of fast manifolds associated with a singularly
perturbed vector field $\Phi(z,\delta)$ the property of minimality
can be written in a rather simple way. If $M_{f}$ is a minimal fast
linear subspace associated with a singularly perturbed vector field
$\Phi(z,\delta)$ then dimension $n_{f}=\dim M_{f}$ cannot be reduced.

Call this minimal subspace $M_{f}$ a linear subspace of fast motions
of $\Phi(z,\delta)$.

\subsection{Fast-slow decomposition of Singularly Perturbed Vector Fields.}

Fix an asymptotic fast-slow vector field $\mathbf{F}(z,x,\delta)$.
Suppose $\{M_{z}\}$ is a fast family associated with $F(z,x.\delta)$
and the fast dimension of $\{M_{z}\}$ is $n_{f}$. Then the vector
field $\mathbf{F}(z,x,\delta)$ is a sum of two vector fields $\mathbf{F}_{f}(z,x,\delta):=Pr_{f}\Phi(z,\delta)+Pr_{f}L(z,x,\delta)$
and $\mathbf{F}_{s}(z,x,\delta):=\mathbf{F}(z,x,\delta)-Pr_{f}\mathbf{F}(z,s,\delta)$.
Here $Pr_{f}\Phi(z,\delta)$ is a projection of $\Phi(z,\delta)$
onto the tangent space $TM_{z}$ of the fast manifold $M_{z}$, $L_{f}(z,x,\delta):=Pr_{f}L(z,x,\delta)$
is the restriction of the linear differential operator $L(z,x,\delta)$
on $TM_{z}$ and $L_{s}(z,x,\delta)$ is a similar projection of $\mathbf{F}(z,x,\delta)$
onto the linear subspace $TM^{z}$ of slow motions that is orthogonal or transverse
to $TM_{z}$.

Call an asymptotic fast-slow vector field $\mathbf{F}(z,x,\delta)$
a uniformly asymptotic fast-slow vector field (or simply a uniform
fast-slow vector field) if 
\[
\lim_{\delta\rightarrow0}\sup_{z\in V;\,x\in V_s}\left|Pr_{s}\mathbf{F}(z,x,\delta)\right|=0.
\]

Denote $\varepsilon:=\sup_{z\in V}\left|Pr_{s}\Phi(z,\delta)\right|$,
which is a new small parameter, $\varepsilon<\varepsilon_{0}=\sup_{z\in V}\left|Pr_{s}\Phi(z,\delta_{0})\right|$;
$F(z,\delta):=Pr_{f}\Phi(z,\delta)$ is the fast sub-field and $G(z,\delta):=\frac{Pr_{s}\Phi(z,\delta)}{\sup_{z\in V}\left|Pr_{s}\Phi(z,\delta)\right|}$
is a slow sub-field of $\Phi(z,\delta)$ for homogeneous system of the source term.
Then the vector field $\mathbf{F}(z,x,\delta)$ can be represented as a linear combination of its
fast and slow sub-fields i.e.

\begin{equation}
\mathbf{F}(z,x,\delta)=\mathbf{F}_{f}(z,x,\delta)+\varepsilon\mathbf{F}_{s}(z,x,\delta).
\label{MainSPVF}
\end{equation}

where

\begin{equation}
\begin{array}{c}
\mathbf{F}_{f}(z,x,\delta)=\left[F(z,\delta)+L_{f}(z,x,\delta)\right],\\
\mathbf{F}_{s}(z,x,\delta)=\left[G(z,\delta)+L_{s}(z,x,\delta)\right].
\end{array}
\label{MainSPVF1}
\end{equation}

Remark the small parameter $\epsilon$ is a function of the small
parameter $\delta$. If $\delta\rightarrow0$, then $\epsilon\rightarrow0$.

For typical reaction-kinetic of combustion systems the situation with Eq. (\ref{MainSPVF})
can be essentially simplified and the regular theory of singularly perturbed
system of ODE can be adapted under the following main assumptions.

\emph{
Main assumptions made:
}
\begin{enumerate}
\item The fast linear operator $L_{f}(z,x,\delta)$ does not depend on
 the small parameter $\delta,$ i.e. to the leading order
it can be written as $L_{f}(z,x)$;
\item Transport processes for the fast and slow variables have the same order,
because diffusion and convection processes do not depend directly on reaction
processes. It means that the fast operator $L_{f}(z,x)$ can be rewritten
as $L_{f}(z,x):=\varepsilon L_{f}(z,x)$ and
\end{enumerate}
\begin{equation}
L_{f}(z,x) \approx L_{s}(z,x) \approx O(1).
\label{Asum2}
\end{equation}

For any practical implementation of the proposed construction of singularly
perturbed vector fields we have to find a way to determine the fast
manifolds. For the moment this can be achieved for the linear case i.e. for
the case where all fast manifolds are parallel to a fixed linear subspace
$L_{f}$. In the next section we shall discuss the linear case of fast-slow
vector fields in more details.

\section{Fast-slow Vector Fields with Linear Fast Subspace}

For any realistic complex model the small parameter $\delta$ is unknown
and this fact restricts possible applications of the proposed asymptotic
theory. In this section the proposed asymptotic theory is adopted and further
developed for practical problem in the simplest possible case of linear fast
manifolds.

Thus, it is assumed any fast manifold at any point $z$ is parallel to a linear subspace
$M_{f}(z)$ with fixed dimension $n_{f}$. Note that for many applications
an assumption that $M_{f}$ does not depend on $z$ is very natural. For instance,
in the case of chemical kinetics, by using mass action law the chemical source term
is represented as composition of linear operator (given by the system stoichiometric
matrix) and non-linear operator describing the rates of elementary reactions.

\subsection{Fast-Slow decomposition.}

Fix a uniformly asymptotic fast-slow vector field $\mathbf{F}(z,x,\delta)$.
Suppose that the fast subspace $M_{f}$ does not depend on $z$ and
$\dim M_{f}=n_{f}$. The vector field $\mathbf{F}(z,x,\delta)$ is
a sum of two vector fields
$\mathbf{F}_{f}(z,x,\delta):=Pr_{f}\mathbf{F}(z,x,\delta)$
and
$\mathbf{F}_{s}(z,x,\delta):=\mathbf{F}(z,x,\delta)-Pr_{f}\mathbf{F}(z,x,\delta)$. 

The uniformity condition permits us to represent a uniformly singularly
perturbed vector field and a corresponding dynamical system (\ref{eq:main})
as a standard singularly perturbed system (SPS) as in the following.

Suppose $u:=Pr_{f}z$ and $v:=Pr_{s}z$ are fast
and slow variables that represent a new coordinate system with $n_{f}$
fast variables $u$ and $n_{s}=n-n_{f}$ slow variables $v$; $\varepsilon:=\sup_{z\in V}\left|Pr_{s}\Phi(z,\delta)\right|$
is a small parameter, $\varepsilon<\varepsilon_{0}=\sup_{z\in V}\left|Pr_{s}\Phi(z,\delta_{0})\right|$;
$F(u,v,\delta)$ is a representation of $F(z,\delta):=Pr_{f}\Phi(z,\delta)$
in the new coordinate system $(u,v)$ and $G(u,v,\delta)$ is a representation
of $G(z,\delta):=\frac{Pr_{s}\Phi(z,\delta)}{\sup_{z\in V}\left|Pr_{s}\Phi(z,\delta)\right|}$
in the new coordinate system $(u,v)$.

Hence the system (\ref{eq:main}) has the standard SPS form 
\begin{equation}
\frac{\partial u}{\partial\tau}=F(u,v,\delta)+L_{f}(u,v;x,\delta)\label{eq:5}
\end{equation}

\begin{equation}
\frac{\partial v}{d\tau}=\varepsilon G(u,v,\delta)+\varepsilon L_{s}(u,v;x,\delta).\label{eq:6}
\end{equation}
 in the new coordinate system $(u,v)$.

Remind once again that the small parameter $\epsilon$ is a function of the
small parameter $\delta$ \cite{BGG2006} and if $\delta\rightarrow0$, then $\epsilon\rightarrow0$.

By re-scaling the time $t=\frac{\tau}{\varepsilon}$ from slow to fast one we
can rewrite the previous fast-slow system in an another equivalent standard form
\begin{equation}
\varepsilon\frac{\partial u}{\partial t}=F(u,v,\delta)+L_{f}(u,v;x,\delta)\label{eq:5-1}
\end{equation}

\begin{equation}
\frac{\partial v}{d\tau}=G(u,v,\delta)+L_{s}(u,v;x,\delta).\label{eq:6-1}
\end{equation}
By using a formal substitution $\varepsilon=0$ we can write an analog
of slow invariant manifold 
\begin{equation}
F(u,v,0)+L_{f}(u,v;x,0)=0.\label{eq:formalIM}
\end{equation}

Under \emph{the main assumptions} above the system Eqs. (\ref{eq:5-1})-(\ref{eq:6-1}) can
be further simplified to

\begin{equation}
\varepsilon\frac{\partial u}{\partial t}=F(u,v,\delta)+\varepsilon L_{f}(u,v,x)\label{eq:5-1-1}
\end{equation}

\begin{equation}
\frac{\partial v}{d\tau}=G(u,v,\delta)+L_{s}(u,v,x).\label{eq:6-1-1}
\end{equation}

Call this system as a reaction-diffusion fast-slow vector field.

By using a formal substitution $\varepsilon=0$ we can write the same
zero approximation (same as for homogeneous sub-system) of the slow invariant manifold,
namely 
\begin{equation}
F(u,v,0)=0.\label{eq:formalIM-1}
\end{equation}

\subsection{Singularly Perturbed Vector Fields: non asymptotic definition.}

In the previous definition of an asymptotic
fast-slow vector field $\mathbf{F}(z,x,\delta)$, a small parameter $\delta$ is unknown.
Meanwhile the main geometrical idea is still useful if some previous knowledge about
a scaling is known. It means that some "small" number $\varepsilon_{0}$
is fixed for corresponding processes (models) and any parameter $\varepsilon<\varepsilon_{0}$
can be considered as a small system parameter.

Suppose a smooth fast-slow vector field $\mathbf{F}(z,x)$ is defined
in a structured domain $V\subset R^{n}$, $z\in V$, in a parametric
domain $V_{s}$, $x\in V_{s}$ and $M_{f}$ is the fast sub-field of Eq. (\ref{MainSPVF}).

Moreover $\sup_{z\in V}\left|Pr_{s}\mathbf{F}(z,x)\right|<\varepsilon_{0}$.

Suppose as well, as in the previous subsection that $u:=Pr_{f}z$
and $v:=Pr_{s}z$ are fast and slow variables that represent a new
coordinate system with $n_{f}$ fast variables $u$ and $n_{s}=n-n_{f}$
slow variables $v$; $\varepsilon:=\sup_{z\in V}\left|Pr_{s}\mathbf{F}(z)\right|$
is a small system parameter; $F(u,v)+L_{f}(u,v,x)$ is a representation
of $Pr_{f}\mathbf{F}(z)$ and $G(u,v)+L_{s}(u,v,x)$ is a representation
of $\frac{Pr_{s}\Phi(z)}{sup_{z\in V}\left|Pr_{s}\Phi(z)\right|}$
in the new coordinate system $(u,v)$.

Hence the system (\ref{eq:main}) can be cast in the following form:

\begin{equation}
\frac{\partial u}{\partial\tau}=F(u,v)+L_{f}(u,v,x),
\label{eq:7}
\end{equation}

\begin{equation}
\frac{\partial v}{d\tau}=\varepsilon\left[G(u,v)+L_{s}(u,v,x)\right],
\label{eq:8}
\end{equation}
in the new coordinate system $(u,v)$ as SPS. In similar manner modifications can be used for
system representations (\ref{eq:5-1})-(\ref{eq:6-1}) and (\ref{eq:5-1-1})-(\ref{eq:6-1-1})
can be defined.

\section{Fast motion time estimates}

In this section a formal definition of slow manifold is justified by estimating
influence of transport for the system (\ref{eq:5-1-1})-(\ref{eq:6-1-1})
that represents a a reaction-diffusion fast-slow vector field.

\subsection{Fast motion time estimates for ODE}

Consider first the system of ordinary differential equations in the standard SPS form

\begin{equation}
\begin{array}{c}
\frac{dx}{dt}=f(x,y),\\
\varepsilon\frac{dy}{dt}=g(x,y).
\end{array}
\end{equation}

where $x\in R^{n},y\in R^{m}.$Here $x$ is a slow vector, $y$ is
a fast vector. Suppose that $(x_{0},y_{0)}$ is an initial data for
this system and $g(x_{0},y_{0})\neq0$ . The subspace $L_{x_{0}}=[(x,y)\in R^{n+m}:x=x_{0}]$
is a fast subspace that contains $(x_{0},y_{0})$. Our main assumption
here is simplicity of the slow invariant manifold $g(x,y)=0$. It
means that the equation $g(x,y)=0$ a zero approximation $\varepsilon=0$
of a stable invariant slow manifold. It means that any fast subspace
has a one point intersection $(x_{0},y_{s})$ with the slow invariant
manifold that is an attractive singular point of the fast sub-system
$\varepsilon\frac{dy}{dt}=g(x,y)$ \cite{BGG2006, BGM2008}. 

Our next assumption used simplicity of fast dynamics. Namely,
a length of the fast trajectory, that joints points $(x_{0},y_{0})$
and the fast singular point $(x_{0},y_{s})$ is less than $2|y_{0}-y_{s}|$.

For any $\varepsilon>0$ introduce the open set $F_{\sqrt{\varepsilon}}:=\{(x,y)\in R^{n+m}|g(x,y)|<\sqrt{\varepsilon}$.
Outside of the slow neighborhood $F_{\sqrt{\varepsilon}}$ of the
slow manifold $F_{s}:=\{(x,y)|g(x,y)=0\}$ the fast component of
the vector field $\Phi(x,y):=\{f(x,y),$$\frac{1}{\varepsilon}g(x,y)$
satisfies to the inequality $\frac{1}{\varepsilon}|g(x,y)\geq\frac{1}{\sqrt{\varepsilon}}$. 

The fast trajectory with the initial point $y_{0}$ is a curve $\varphi:[0,\infty)\to L_{x_{0}}$where
$0\leq t<\infty$ and $\varphi'(t)=\frac{1}{\varepsilon}g(x_{0},\varphi(t))$.
Under our assumptions its length 
\[
l_{\varphi}:=\int_{0}^{\infty}|\varphi'(t)|dt=\int_{0}^{\infty}\frac{1}{\varepsilon}|g(x_{0},\varphi(t))|dt\leq2|y_{0}-y_{s}|.
\]

For any $t_{0}>\sqrt{\varepsilon}2|y_{0}-y_{s}|$ we have 
\[
l_{\varphi}=\int_{0}^{t_{0}}\frac{1}{\varepsilon}|g(x_{0},\varphi(t))|dt\geq\int_{0}^{t_{0}}\frac{1}{\sqrt{\varepsilon}}dt=2|y_{0}-y_{s}|.
\]

It means that for any $t_{0}>\sqrt{2\varepsilon}|y_{0}-y_{s}|$ the
point $\varphi(t_{0})$ belongs to the slow neighborhood $F_{\sqrt{\varepsilon}}$.
Therefore the fast motion time is less than $\sqrt{2\varepsilon}2|y_{0}-y_{s}|$.
After this time the solution of the fast subsystem $\varepsilon\frac{dy}{dt}=g(x,y)$
belongs to the slow neighborhood $F_{\sqrt{\varepsilon}}$. The influence
of the slow sub-system to this estimate is negligible.

\subsection{Fast motion time estimates for models with diffusion.}

Consider the system of PDEs with 1D spatial transport/diffusion terms. Under main assumption 
(1), (2) above, it can be cast in the fast time as

\begin{equation}
\begin{array}{c}
\frac{du(x,t)}{dt}=f\left(u(x,t),\,v(x,t)\right)+L_{1,x}(u(x,t),\,v(x,t))\\
\frac{dv(x,t)}{dt}=\frac{1}{\epsilon}g\left(u(x,t),\,v(x,t)\right)+L_{2,x}(u(x,t),\,v(x,t))
\end{array}
\end{equation}

Here $0\leq x\leq1$ and $L_{1,x},L_{2,x}$ are elliptic differential
operators of the second order.

The transport term is treated as slow compared to the fast component
of the vector field, i.e

\[
|L_{1,x}(u(x,t),\,v(x,t))|\leq K|g\left(u(x,t),\,v(x,t)\right)|
\]

\[
|L_{2,x}(u(x,t),\,v(x,t))|\leq K|g\left(u(x,t),\,v(x,t)\right)|
\]
outside of the slow neiboobhood $F_{\sqrt{\varepsilon}}$ . Here $K$
is a constant that typically do not exceed $3$.

An additional assumption for the fast subsytem is 
\[
|L_{2,x}(u(x,0),\,v(x,t))|\leq K|g\left(u(x,t),\,v(x,t)\right)|,
\]
for any $x\in[0,1]$ and any $t\in[0,\infty)$.

The slow system evolution is then controlled by 
\[
u(x,t)=\left(u_{1}(x,t),...,\,u_{m_{s}}(x,t)\right),
\]

which are assumed to change slowly comparatively to the fast variables
\[
v(x,t)=\left(v_{1}(x,t),...,\,v_{m_{f}}(x,t)\right),\quad m_{s}+m_{f}=n.
\]

The transport diffusion terms are represented first by very general
and smooth differential operators $L_{1,x}(u(x,t),\,v(x,t)),L_{2,x}(u(x,t),\,v(x,t))$

Initial data for the system are 
\begin{equation}
u(x,0)=u_{0}(x),v(x,0)=v_{0}(x).
\label{eq:indata}
\end{equation}

Our main goal is to estimate influence of the transport operators to
the fast time estimates obtained in the previous section. 

Fix $x_{0}\in[0,1]$ and check length of the fast trajectory that
belongs to the fast subspace $L_{u(x_{0},0)}$, its starting point
is $y_{0}:=v(x_{0},0)$ and its final point is $y_{s}:=L_{u(x_{0},0)}\cap F_{s}$.

The fast trajectory with the initial point $y_{0}$ is a curve $\varphi:[0,\infty)\to L_{u(x_{0},0)}$ where
$0\leq t<\infty$ and $\varphi'(t)=\frac{1}{\varepsilon}g(x_{0},\varphi(t)+L_{2,x}(u(x,t),\,v(x,t)))$.
Under our assumptions its length 
\[
l_{\varphi}:=\int_{0}^{\infty}|\varphi'(t)|dt=\int_{0}^{\infty}\left|\frac{1}{\varepsilon}g(u(x_{0},0),\varphi(t))+L_{2,x}(u(x_{0},0),\,v(x_{0},t)))\right|dt\leq
\]
\[
\int_{0}^{\infty}\left|\frac{1}{\varepsilon}g(u(x_{0},0),\varphi(t))\right|dt+\int_{0}^{\infty}\left|L_{2,x}(u(x_{0},0),\,v(x_{0},t)))\right|dt\leq2|y_{0}-y_{s}|
\]
\[
+\int_{0}^{\infty}K\left|g(u(x_{0},0),\varphi(t))\right|dt\leq2|y_{0}-y_{s}|+\varepsilon K2|y_{0}-y_{s}|=2(1+\varepsilon K)|y_{0}-y_{s}|.
\]

It means that for any $t_{0}>\sqrt{2\varepsilon}(1+\varepsilon K)|y_{0}-y_{s}|$
the point $\varphi(t_{0})$ belongs to the slow neighborhood $F_{\sqrt{\varepsilon}}$.
Therefore, the fast motion time is less than $\sqrt{2\varepsilon}(1+\varepsilon K)|y_{0}-y_{s}|$.
After this time the solution of the fast subsystem
$\frac{dy}{dt}=\frac{1}{\epsilon}g\left(u(x_{0},0),\,v(x_{0},t)\right)+L_{2,x}(u(x_{0},0),\,v(x_{0},t))$
belongs to the slow neighborhood $F_{\sqrt{\varepsilon}}$. The influence
of the slow sub-system to this estimate is negligible similar as in the previous subsection.


\section{Singularly perturbed profiles and the REDIM approach}

In this section the REDIM method is discussed as a method to
construct the manifold approximating relatively slow evolution of the detailed
system solution profiles. Recall definition of singularly perturbed
profiles \cite{BCMGM2016}. Accordingly, the following
representation of the system Eq. (\ref{eq:main}) can be obtained

\begin{equation}
\left\{ \begin{array}{c}
\frac{du(x,t)}{dt}=F_{s}\left(u(x,t),\,v(x,t)\right)+L_{1,x}(u(x,t),\,v(x,t))\\
\frac{dv(x,t)}{dt}=\frac{1}{\epsilon}F_{f}\left(u(x,t),\,v(x,t)\right)+L_{2,x}(u(x,t),\,v(x,t))
\end{array}\right.
\label{eq:SPP1}
\end{equation}

The slow system evolution is then controlled by 
\[
u(x,t)=\left(u_{1}(x,t),...,\,u_{m_{s}}(x,t)\right),
\]

which are assumed to change slowly comparatively to the fast variables
\[
v(x,t)=\left(v_{1}(x,t),...,\,v_{m_{f}}(x,t)\right),\quad m_{s}+m_{f}=n.
\]
 
We suppose that $u(x,t)$, $v(x,t)$ are smooth functions.
Initial data for the system Eq. (\ref{eq:SPP1}) are

\begin{equation}
u(x,0)=u_{0}(x),v(x,0)=v_{0}(x).
\label{eq:indata-1}
\end{equation}

Recall that functions $F_{s},F_{f}$ are of the same order.
Then $\left|\left|\frac{dU}{dt}\right|\right|\sim O\left(1\right)$,
while $\left|\left|\frac{dV}{dt}\right|\right|\sim O\left(\frac{1}{\varepsilon}\right)$.
Suppose also that operators (see the assumption above (2))
$$L_{1,x}(u(x,t),v(x,t)),L_{2,x}(u(x,t),v(x,t))$$ have the same order as $F_{s},F_{f}$ terms.

Recall that the zero approximation $S$ of the slow invariant manifold
in the phase space $(u,v)$ (the space of species) is represented
in the implicit form 
\[
F_{f}(u,v)=0.
\]
The initial profile is $\varGamma_{0}(x):=(u_{0}(x),v_{0}(x));\,u_{0}(x)=u(x,0),\,v_{0}(x)=v(x,0)$.
Denote $\varGamma(x,t)$ a profile that is the solution of (\ref{eq:SPP1})
at time $t$ with the initial profile (initial data) $\varGamma_{0}(x)$.

For a system in the general form (\ref{eq:main}) this information is absent, thus, the question is how to access
\begin{equation}
F_{f}\left(u(x,t),\,v(x,t)\right)=0,
\label{eq:RE-0}
\end{equation}
as e.g. the zero order approximation $\varGamma_{0}(x,t)$ of $\varGamma(x,t)$
which belongs to $S$ for all $t$, represents the main problem of model reduction for
a reaction-diffusion system.

The set $RM:=\cup_{t\in(0,\infty)}\varGamma(x,t)$ is called the reaction-diffusion
manifold (REDIM) and $RM_{0}:=\cup_{t\in(0,\infty)}\varGamma_{0}(x,t)$
is its zero approximation (for $\epsilon=0$).

Note that if the dimension of the profile is equal to $s$ ($\dim\varGamma(x,t)=s$),
then $\dim RM=\dim RM_{0}\leq s+1$.

\subsection{REDIM}

In the framework of the REDIM \cite{BM2007}, the manifold of the relatively slow profile evolution
$RM_{0}$ is constructed / approximated by using the so-called \emph{Invariance condition} (see e.g. \cite{GKZ02,
GKZ04, BM2007} for more details). The construction of an explicit representation of a low-dimensional manifold
\begin{equation}
RM_{0} = \{ z: z=z(\theta), \theta \in R^{m_{s}}\},
\label{eq:Manifold}
\end{equation}
starts from an initial solution $z=z_0(\theta)$ and then it is integrated with the vector field of the PDEs
reaction-diffusion system:
\begin{equation}
\frac{\partial z(\theta)}{\partial \tau}=(I-z_{\theta}z_{\theta}^+)(\Phi(z(\theta),\delta)+L(z(\theta),x,\delta)),
\label{eq:M_REDIM}
\end{equation}
where the evolution of the manifold along its tangential space is forbidden by restricting it to the normal
(or transverse) subspace. This is achieved by the local projector: $Pr_{TM^\perp}=(I-z_{\theta}z_{\theta}^+)$,
here $I$ identity matrix, $z_{\theta}$ denotes the tangential subspace and $z_{\theta}^+$ is the Moore-Penrouse
pseudo-inverse of the local coordinates Jacobi matrix $z_{\theta}$.
In this special case the evolution of the manifold Eq. (\ref{eq:M_REDIM}) is computed in the normal direction until
the stationary solution is reached \cite{BM2007}.

Now, if the main assumption of the study is valid, the manifold will evolve within fast manifolds of the vector
field Eq. (\ref{eq:main}) and will converge asymptotically to an invariant system manifold $RM_{0}$ approximating the slow profile
evolution \cite{BCMGM2016}.


\begin{figure}[ht]
\centering
\includegraphics[scale=0.6]{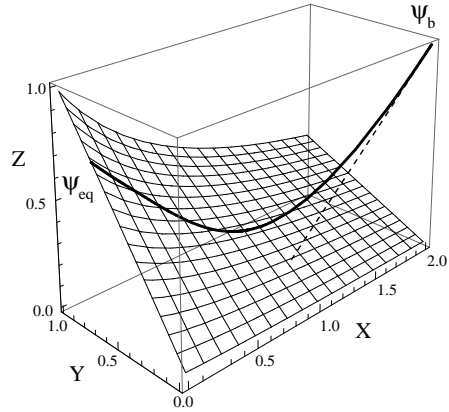}
\caption{System state space (X,Y,Z) is shown. 2D slow manifold for the pure homogeneous system
(\ref{eq:M_MentenODE}) is represented by a mesh. System stationary solution profile (\ref{eq:M_MentenPDE}),
black solid curve) in 1D case corresponds to 1D REDIM due to dimensional 
considerations. The approximation of the fast part of the homogeneous system (\ref{eq:M_MentenODE})
solution trajectory starting from the boundary state (\ref{eq: for 19}) is shown by the dashed line.}
\label{fig:F1}
\end{figure}

\section{Analysis of 3D Michaelis-Menten model with the Laplacian Operator}

The 3D Michaelis-Menten model is considered here as illustrative example of the REDIM approach.
The original mathematical model of the enzyme biochemical system consists of three ODEs

\begin{align}
 & \frac{dX}{dt}=-XZ+L_{1}(1-Z-\mu(1-Y))\\
 & \frac{dY}{dt}=-L_{3}YZ+\frac{L_{4}}{L_{2}}(1-Y)\\
 & \frac{dZ}{dt}=\frac{1}{L_{2}}((-XZ+1-Z-\mu(1-Y))+\mu)-L_{3}YZ+\frac{L_{4}}{L_{2}}(1-Y)))
\label{eq:M_MentenODE}
\end{align}

The system parameters are taken as $L_{1}=0.99$, $L_{2}=1$, $L_{3}=0.05$,
$L_{4}=0.1$, $\mu=1$ (see e.g. \cite{BCMGM2016, RF01} for details and references).
By taking the 1D diffusion into account we obtain the following PDEs
system with the constant diffusion coefficient was taken as $\delta=0.01$:

\begin{align}
 & \frac{\partial X}{\partial t}=-XZ+L_{1}(1-Z-\mu(1-Y))+\delta\Delta X\\
 & \frac{\partial Y}{\partial t}=-L_{3}YZ+\frac{L_{4}}{L_{2}}(1-Y)+\delta\Delta Y\\
 & \frac{\partial Z}{\partial t}=\frac{1}{L_{2}}((-XZ+1-Z-\mu(1-Y))+\mu)-L_{3}YZ+\frac{L_{4}}{L_{2}}(1-Y)))+\delta\Delta Z
\label{eq:M_MentenPDE}
\end{align}

The system (\ref{eq:M_MentenPDE}) is considered
with the following initial and boundary conditions:

\begin{align}
 & \begin{pmatrix} 
 X(t,0)=X_{eq} \\
 Y(t,0)=Y_{eq} \\
 Z(t,0)=Z_{eq}\\
 \end{pmatrix}\label{eq: for 18}\\
 & \begin{pmatrix} 
 X(t,1)=2 \\
 Y(t,1)=0 \\
 Z(t,1)=1\\
 \end{pmatrix}\label{eq: for 19}\\
 & \begin{pmatrix} 
 X(0,x)=(2-X_{eq})x+X_{eq} \\
 Y(0,x)=(-Y_{eq})x+Y_{eq} \\
 Z(0,x)=(1-Z_{eq})x+Z_{eq}\\
 \end{pmatrix}\label{for 20}
\end{align}

Here $(X_{eq},Y_{eq},Z_{eq})$ are coordinates of the equilibrium point and $x$ is spatial variable.
Initial conditions are chosen to be a straight lines, they satisfy the general assumption - join initial
and equilibrium values on the boundaries. 

First, several numerical experiments were performed (see Fig. \ref{fig:F1}). A 2D slow manifold for homogeneous system
(\ref{eq:M_MentenODE}) was found by Global Quasi-Linearisation (GQL) method \cite{BGG2006} 
(see Appendix for a short description of GQL). Stationary system (\ref{eq:M_MentenPDE}) solution profile was
also integrated. Figure \ref{fig:F1} shows a connection between
the zero approximation of the slow manifold and the profile of the stationary system
solution of the PDE in the original coordinates $(X,Y,Z)$. In Fig. \ref{fig:F1}
the system solution profile can be roughly subdivided into two parts: the slow part
of the stationary solution that is very close to the slow manifold of the homogeneous system
and second one, which is influenced by the diffusion term. The dashed line in this figure
represents an approximation of linear fast sub-field (1D in this case).

\begin{figure}[ht]
 \centering
 \includegraphics[scale=0.6]{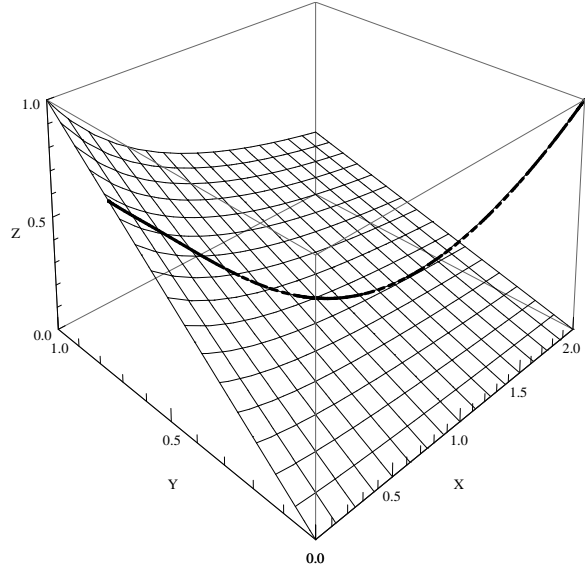}
\caption{REDIM manifold (dashed line), exact stationary solution of the original PDE system (think line).}
\label{fig:F2}
\end{figure}

As in the previous section the main assumption remains the transport term is slow compared
with the fast vector field. By applying the REDIM approach the stationary solution of the
following system should represent the one-dimensional REDIM.
\begin{align}
 & \frac{\partial \Psi}{\partial t}=(I-\Psi_{\theta}\Psi_{\theta}^+)(F_R+F_D)\label{eq:M_MentenREDIM}\\
\end{align}
where following notations have been used, $I$ is 3x3 identity matrix, the system state vector
\[
\Psi=
  \begin{pmatrix}
    X  \\
    Y  \\
    Z
  \end{pmatrix}
 ,\]
and projection matrix to the manifolds' tangent space is given by
\[
\Psi_\theta\Psi_\theta^+=\frac{1}{X_\theta^2+Y_\theta^2+Z_\theta^2}
  \begin{pmatrix}
    X_\theta^2 & X_\theta Y_\theta &  X_\theta Z_\theta  \\
    Y_\theta X_\theta & Y_\theta^2 &  Y_\theta Z_\theta   \\
     Z_\theta X_\theta &  Z_\theta Y_\theta &  Z_\theta^2
  \end{pmatrix}
, \]
and vector fields of reaction and diffusion terms
 \[
L\left(\Psi(\theta)\right)=\delta \theta_x^2
  \begin{pmatrix}
    X_{\theta \theta}  \\
    Y_{\theta \theta}  \\
    Z_{\theta \theta}
  \end{pmatrix}
 ,\]
\[
\Phi\left(\Psi(\theta)\right)=
  \begin{pmatrix}
   -XZ+L_{1}(1-Z-\mu(1-Y)) \\
 -L_{3}YZ+\frac{L_{4}}{L_{2}}(1-Y) \\
 \frac{1}{L_{2}}((-XZ+1-Z-\mu(1-Y))+\mu)-L_{3}YZ+\frac{L_{4}}{L_{2}}(1-Y)))
  \end{pmatrix}
 .\]
Here $\theta$ is the manifold parameter and $\theta_x$ is the gradient of the manifold parameter in $L\left(\Psi(\theta)\right)$. 
Now by using $\theta=X$ as a local manifold parameter, the system (\ref{eq:M_MentenREDIM}) can be simplified
to only two equations for $Y=Y(\theta)$ and for $Z=Z(\theta)$.
They were integrated and the stationary solution has been found for 1D REDIM, which is completely coincides
with the system stationary profile (see Fig. \ref{fig:F2}).

The stationary solution of the system (\ref{eq:M_MentenREDIM}) whith $\theta=(\theta_1,\theta_2)$ represent 2D REDIM.
Here $(\theta_1$,$\theta_2)$ are two manifold parameters. In this case projection matrix of the manifold tangent space 
is considered by $\Psi_\theta \Psi_\theta^+$ where
\[
\Psi_\theta=
  \begin{pmatrix}
    X_{\theta_1} & X_{\theta_2}  \\
    Y_{\theta_1} & Y_{\theta_2}   \\
    Z_{\theta_1} & Z_{\theta_2}
  \end{pmatrix}
, \]
and $\Psi_{\theta}^{+}$ is the Moore-Penrose pseudo-inverse of $\Psi_\theta$:
$ \Psi_{\theta}^{+} =(\Psi_{\theta}^T \Psi_{\theta})^{-1}\Psi_{\theta}^T$.

The components of the diffusion term $L\left(\Psi(\theta)\right)$ are the following:
\[
L\left(\Psi(\theta)\right)^X=\delta(\theta_{1x}, \theta_{2x})
       \begin{pmatrix}
         X_{\theta_1 \theta_1} & X_{\theta_1\theta_2}  \\
         X_{\theta_2 \theta_1} & X_{\theta_2 \theta_2}   \\
            \end{pmatrix}
                \begin{pmatrix}
                        \theta_{1x}  \\
                         \theta_{2x}   \\
                                   \end{pmatrix}
 ,   \]  
 
\[
L\left(\Psi(\theta)\right)^Y=\delta(\theta_{1x}, \theta_{2x})
       \begin{pmatrix}
         Y_{\theta_1 \theta_1} & Y_{\theta_1\theta_2}  \\
        Y_{\theta_2 \theta_1} & Y_{\theta_2 \theta_2}   \\
            \end{pmatrix}
                \begin{pmatrix}
                        \theta_{1x}  \\
                         \theta_{2x}   \\
                                   \end{pmatrix}
 ,   \]    

\[
L\left(\Psi(\theta)\right)^Z=\delta(\theta_{1x}, \theta_{2x})
       \begin{pmatrix}
         Z_{\theta_1 \theta_1} & Z_{\theta_1\theta_2}  \\
         Z_{\theta_2 \theta_1} & Z_{\theta_2 \theta_2}   \\
            \end{pmatrix}
                \begin{pmatrix}
                        \theta_{1x}  \\
                         \theta_{2x}   \\
                                   \end{pmatrix}
 .   \]
 
\begin{figure}[ht]
\centering
\includegraphics[width=0.3\textwidth]{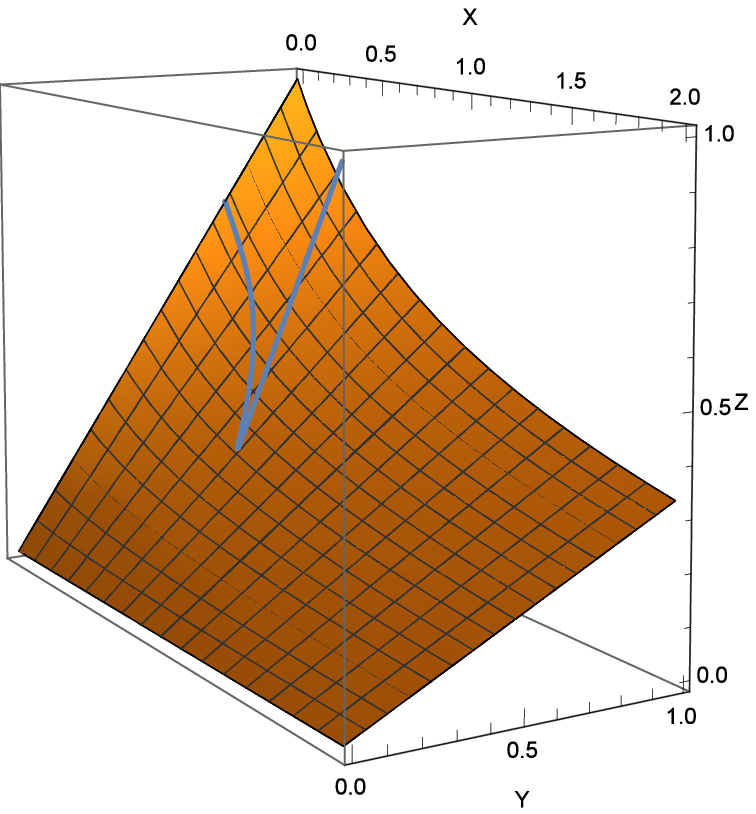}
\includegraphics[width=0.3\textwidth]{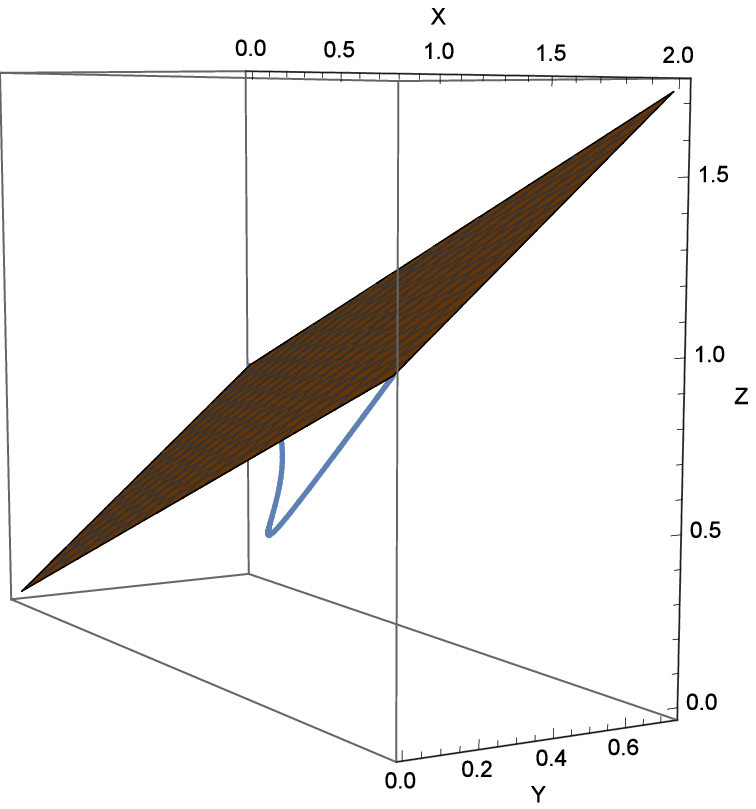}
\includegraphics[width=0.3\textwidth]{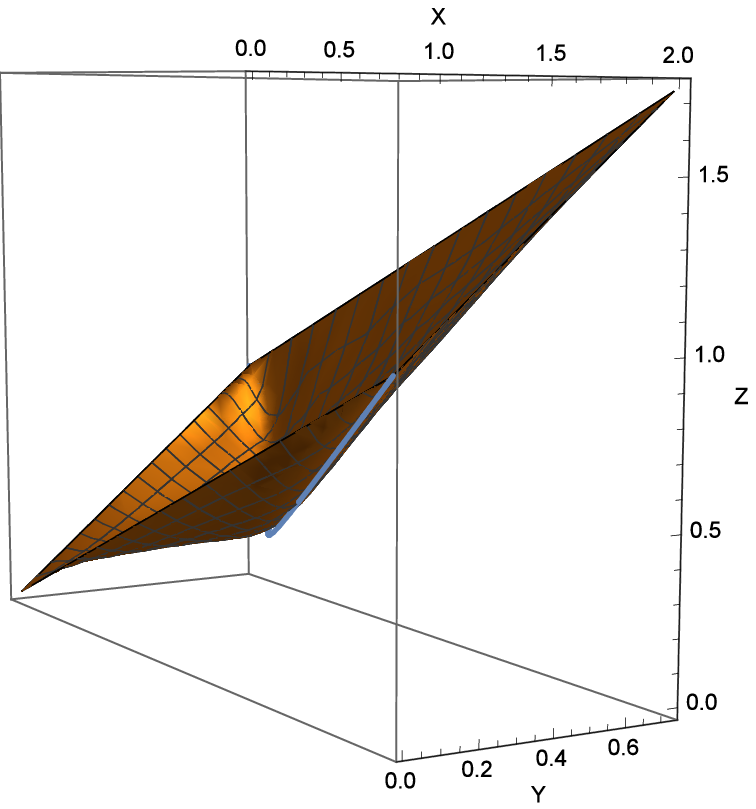}
\caption{On the left: 2D slow homogeneous system manifold, in the middle: an initial solution for the REDIM and
2D REDIM manifold (shown on the right), exact stationary solution of the original PDE system (shown by thick line).}
\label{fig:F3}
\end{figure}

By using $\theta_1=X$, $\theta_2=Y$ as a local coordinates on the manifold of the system (\ref{eq:M_MentenREDIM})
can be simplified to only one equation for $Z=Z(\theta_1,\theta_2)$. It was integrated and the stationary solution
has been found for 2D REDIM.    
Figure \ref{fig:F3} shows a connection between the 2D REDIM, initial solution for the REDIM and/or slow
homogeneous system manifold as in Figs.\ref{fig:F1} and \ref{fig:F2}. The stationary solution profile of
the system illustrates the implementation and quality of the the REDIM approach to approximate the low-
dimensional invariant manifold of relatively slow evolution of the reacting-diffusion system.

Figure \ref{fig:F3} shows the stages of the REDIM construction. On the left the zero order approximation for
a homogeneous system Eq. (\ref{eq:M_MentenODE}). In the middle one can see the stationary solution profile of the
PDEs system Eq. (\ref{eq:M_MentenPDE}), and on the right the converged stationary REDIM equation Eq. (\ref{eq:M_MentenREDIM})
solution is shown together with the stationary systems solution profile. One can see that 2D REDIM manifold approximates
the relatively slow 2D system profile evolution. It means that when the system solution profile evolves Eq. (\ref{eq:M_MentenREDIM})
far form this surface it will evolve relatively fast (see subsection 4.1) towards 2D REDIM along the fast direction of
the fast subspace (see Fig. \ref{fig:F1}) and then finally attains the stationary system solution profile.
In this way, relative fast system dynamics cab be decoupled and the model is reduced to 2D model as a profile evolving
within 2D REDIM.
 
 
\section{Appendix(GQL and system decomposition)} 
Fast sub-fields and fast manifolds play a pivotal role in the theory and applications
 of the SPVF.  The fast manifolds' approximation is
crucial for practical realization of the suggested SPVFs framework. A procedure
for evaluation of the dimension and structure of fast sub-fields is proposed in this
section.
 
In the case when fast manifolds and the system decomposition have
linear structure they can be identified by a gap between the eigenvalues
of an appropriate global linear approximation of the Right Hand Side
(RHS) - vector function of a homogeneous system $\frac{d\psi}{dt}=F(\psi)$ (see \cite{BM2009a} for
detailed discussion)

\[
T\psi\approx F\left(\psi\right).
\]

Note that we did not use a hidden small parameter $\delta$ in $F\left(\psi\right)$,
because its existence is not known 'a priori' and has to be validated
in a course of application of the GQL. Now, if $T$ has two groups
of eigenvalues: so-called small eigenvalues $\lambda\left(\Lambda_{s}\right)$
and large eigenvalues $\lambda\left(\Lambda_{f}\right)$  that have sufficiently
different order of magnitude, then the vector field $F(\psi)$ is
regarded as linearly decomposed asymptotic singularly perturbed vector
field \cite{BGG2006}. Accordingly, fast and slow invariant sub-spaces
given by columns of the matrices $Z_{f},\:Z_{s}$ corresponding 
\cite{MP92} define the slow and  variables . Namely, 

\begin{equation}
T\equiv\begin{pmatrix}Z_{f} & Z_{s}\end{pmatrix}\cdot\begin{pmatrix}\Lambda_{f} & 0\\
0 & \Lambda_{s}
\end{pmatrix}\cdot\begin{pmatrix}\tilde{Z_{f}}\\
\tilde{Z_{s}}
\end{pmatrix},\label{eq:11_T_GQL}
\end{equation}

now, if we denote 

\[
\tilde{Z}=Z^{-1}=\begin{pmatrix}Z_{f} & Z_{s}\end{pmatrix}^{-1}=\begin{pmatrix}\tilde{\left(Z_{f}\right)}_{m_{s}\times n}\\
\left(\tilde{Z_{s}}\right)_{m_{f}\times n}
\end{pmatrix},
\]

then, new coordinates suitable for an explicit decomposition (and
coordinates transformation) are given by $(U,\,V)$:

\begin{equation}
\begin{array}{c}
U:=\tilde{Z}_{f}\,\psi\\
V:=\tilde{Z}_{s}\,\psi
\end{array}.\label{eq:13_VARIABLES-2}
\end{equation}

The decomposed form and corresponding fast and slow subsystems becomes

\begin{equation}
\left\{ \begin{array}{c}
\frac{dU}{dt}=\tilde{Z}_{f}\cdot F\left(\left(Z_{f}\:Z_{s}\right)\left(\begin{array}{c}
U\\
V
\end{array}\right)\right)\\
\frac{dV}{dt}=\tilde{Z}_{s}\cdot F\left(\left(Z_{f}\:Z_{s}\right)\left(\begin{array}{c}
U\\
V
\end{array}\right)\right)
\end{array}\right..\label{eq:14_ODE_DECOMPOSED}
\end{equation}

The small system parameter controlling the characteristic time scales
in (\ref{eq:14_ODE_DECOMPOSED}) can be estimated by the gap between
the smallest eigenvalue of the slow  group and the largest eigenvalue
of the fast group of eigenvalues \cite{BM2009}

\begin{equation}
\varepsilon=\frac{max\left|\lambda\left(\Lambda_{s}\right)\right|}{min\left|\lambda\left(\Lambda_{f}\right)\right|}\ll1.\label{eq:15_SMALL}
\end{equation}

In principle, the idea of the linear transformation is not new, see
e.g. \cite{OM1998}, but the principal point of the developed algorithm
concerns evaluation of this transformation. We have developed the
efficient and robust method that produces the best possible (to the
leading order) decomposition with respect to existing multiple-scales
hierarchy (see the attachment and \cite{BG2013,BGM2008,BM2009} for
more details).
 
\section{Conclusions}
The framework for manifolds based model reduction of the reaction-diffusion system has been established in
the current work. This follows the original ideas of the singularly perturbed vector fields developed earlier.
Within the suggested concept the problem of model reduction is treated as restriction of the original system
to a low-dimensional manifold embedded in the systems state space. The manifold encounters the stationary states
of the degenerate fast sub-field of the vector field defined by the reaction-diffusion system.

The main assumption of weak dependence of the fast system sub-filed of the reaction-diffusion PDEs vector field
on the diffusion has been formulated.
Under this assumption the theory of singularly perturbed vector fields was extended to the the systems with the
molecular diffusion included. The developed framework can be used to justify the so-called REDIM method developed
for reacting flow systems. For illustration Michaelis-Menten chemical kinetics model is extended to describe
reaction-diffusion process. This example is used as an application that illustrate the method and the suggested
framework. It was found that relatively fast 1D sub-field can be decoupled and the system can be reduced and
represented by 2D reduced system.

\section*{Acknowledgments} Financial support by the DFG within the German-Israeli 
Foundation under Grant GIF (No: 1162-148.6/2011) is gratefully acknowledged.

\end{document}